\documentclass[12pt]{amsart}
\usepackage{epsfig}
\usepackage{graphics}
\numberwithin{equation}{section}
\theoremstyle{plain}
   \newtheorem{thm}{Theorem}[section]
   \newtheorem{cor}[thm]{Corollary}
   
   \newtheorem{lemma}[thm]{Lemma}

   \newtheorem{q}[thm]{Question}
   
\theoremstyle{definition}
   \newtheorem{dfn}[thm]{Definition}
   
\theoremstyle{remark}
   \newtheorem{rem}[thm]{Remark}

\newcommand{\R}{\mathbb R}

\newcommand{\ve}{\varepsilon}
\newcommand{\bd}{\partial}

\begin{document}
\title{Legendrian graphs and quasipositive diagrams}
\author{Sebastian Baader and Masaharu Ishikawa}
\footnote[0]{MI supported by MEXT, Grant-in-Aid for Young Scientists (B) (No. 16740031).}
\address{Department of Mathematics, ETH Z\"urich, Switzerland}
\email{sebastian.baader@math.ethz.ch}
\address{Department of Mathematics, Tokyo Institute of Technology,
2-12-1, Oh-okayama, Meguro-ku, Tokyo, 152-8551, Japan}
\email{ishikawa@math.titech.ac.jp}

\begin{abstract}
In this paper we clarify the relationship between ribbon surfaces of
Legendrian graphs and quasipositive diagrams by using
certain fence diagrams.
As an application, we give an alternative proof of a theorem
concerning a relationship between quasipositive fiber surfaces
and contact structures on $S^3$. We also answer a question of L.~Rudolph
concerning moves of quasipositive diagrams.
\end{abstract}

\maketitle


\section{Introduction}

\begin{figure}[htbp]
   \centerline{\input{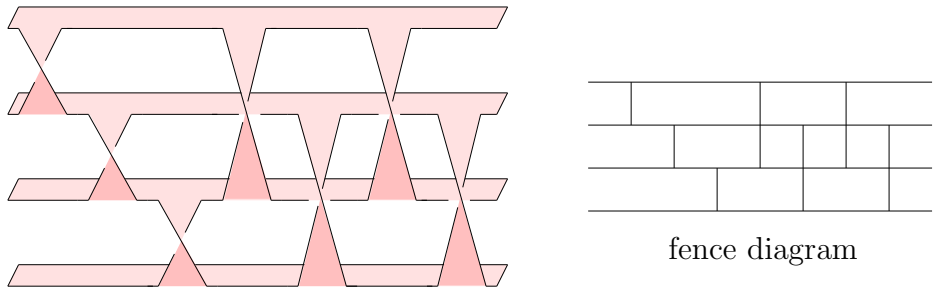}}
   \caption{An example of quasipositive surface. The boundary 
is the knot $10_{145}$ in Rolfsen's notation~\cite{rolfsen}.\label{fig7}}
\end{figure}

A link is called {\it quasipositive} if it has a diagram
which is the closure of a product of conjugates
of the positive generators of the braid group.
If the product consists only of words of the form
\[
   \sigma_{i,j}=(\sigma_i\cdots\sigma_{j-2})\sigma_{j-1}
   (\sigma_i\cdots\sigma_{j-2})^{-1}
\]
then the link obtained as its closure is called
{\it strongly quasipositive}.
Let $b$ be the braid index of a braid diagram of a quasipositive link.
The link spans a canonical Seifert surface
consisting of $b$ copies of disjoint parallel disks with
a band for each $\sigma_{i,j}$, for example see the diagram on the left
in Figure~\ref{fig7}.
We call such a diagram of a Seifert surface a {\it quasipositive diagram}
and each band a {\it positive band}.
We say a Seifert surface is {\it quasipositive} if
it has a quasipositive diagram.
On the right in Figure~\ref{fig7},
the quasipositive diagram is represented by using a graph,
which is called a {\it fence diagram}.

\begin{figure}[htbp]
   \centerline{\input{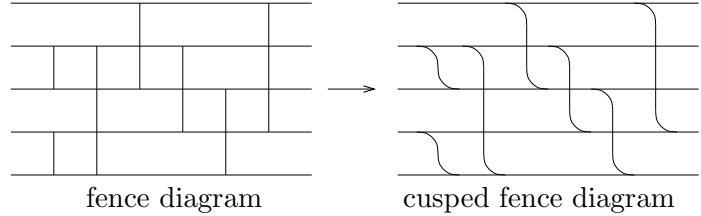}}
   \caption{From a fence diagram to a cusped fence diagram,
which is a front projection of a Legendrian graph.\label{figure22}}
\end{figure}

To relate fence diagrams to contact topology,
we replace each endpoint of vertical lines
as shown in Figure~\ref{figure22}
and call the obtained diagram its {\it cusped fence diagram}.
A cusped fence diagram is regarded as a front projection
of a Legendrian graph. We will see
that the Legendrian ribbon of this Legendrian graph is
the quasipositive surface of the fence diagram (Lemma~\ref{lemma2}).
Conversely, when a front projection of a Legendrian graph is given,
we can make a fence diagram whose quasipositive surface
is a Legendrian ribbon of the given Legendrian graph.
In particular, the Legendrian ribbon of a Legendrian graph 
in $\R^3$ with the standard contact structure is
a quasipositive surface (Theorem~\ref{prop1}).

As an application, we give an alternative proof
of~\cite[Proposition~2.1 (1)$\Leftrightarrow$(2)]{hedden} which
states that a fiber surface is quasipositive if and only if
it supports the standard contact structure on $S^3$
(Theorem~\ref{thmhedden}).
The proof in~\cite{hedden} is based on the work of
E.~Giroux~\cite{giroux} and L.~Rudolph~\cite{rudolph:98},
and uses plumbings. Hence it requires the affirmative answer to
J.~Harer's conjecture~\cite{harer:82}
due to Giroux and N.~Goodman~\cite{giroux, goodman}.
In our proof, one direction follows from the fence diagram argument
and the other one is done according to an argument in~\cite{etnyre}.
In particular, both directions do not need plumbings.

\begin{figure}[htbp]
   \centerline{\input{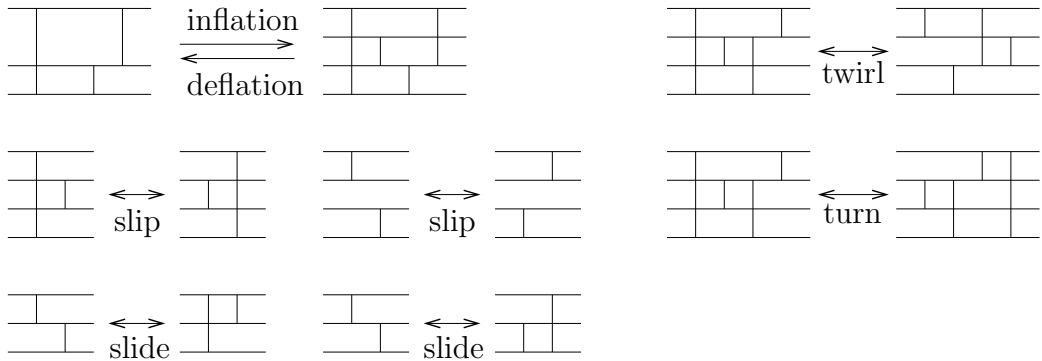}}
   \caption{Fundamental moves of quasipositive diagrams.
Each positive band may pass over several horizontal lines.\label{figure3}}
\end{figure}

Next we observe Legendrian isotopy moves of Legendrian graphs
using moves of quasipositive surfaces.
In~\cite{rudolph:98}, Rudolph introduced fundamental moves
of quasipositive diagrams
named {\it inflations, deflations, slips, slides, twirls and turns}.
These moves are summarized in Figure~\ref{figure3}.
The same figures can be found in his paper~\cite{rudolph:98}
with more precise definitions.
In his notation, we only consider the case where the signs
$\ve$ assigned to bands are positive.
There is a remark about these definitions, see Remark~\ref{remdfn} below.
We will prove that all Legendrian isotopy moves of
Legendrian graphs are expressed by 
inflations, deflations, slips and slides
of fence diagrams (Theorem~\ref{thmiso2}).

In the last two sections we study quasipositive annuli.
We first prove that the moves of fence diagrams
of quasipositive annuli correspond to Legendrian isotopy moves
(Theorem~\ref{thmiso}).
It is important to remark that
the same assertion is not true for quasipositive surfaces.
Secondly, we study
the Thurston-Bennequin invariant and the rotation number
of fence diagrams of quasipositive annuli, which are
defined by those of cusped fence diagrams.
We prove that the Thurston-Bennequin invariant and the rotation number 
of a fence diagram of a quasipositive annulus are invariant under
inflations, deflations, slips, slides, twirls and turns
(Theorem~\ref{mainthm}). As a corollary, we conclude that
there exists a quasipositive surface with
two different quasipositive diagrams which are not related by
inflations, deflations, slips, slides, twirls and turns.
This answers a question of Rudolph in~\cite[Remark on p.263]{rudolph:98}.

From the results in this paper, we conclude that
there exist surjective, non-injective maps
\[
\begin{split}
\{\text{trivalent Legendrian ribbons}\}_{/\sim} 
&\to 
\{\text{quasipositive diagrams}\}_{/\sim} \\
&\to
\{\text{quasipositive surfaces}\}_{/\sim},
\end{split}
\]
where $\{$trivalent Legendrian ribbons$\}_{/\sim}$
is the class of trivalent Legendrian graphs
up to Legendrian isotopy,
$\{$quasipositive diagrams$\}_{/\sim}$ is
the class of quasipositive diagrams up to
inflations, deflations, slips, slides, twirls and turns,
and
$\{$quasipositive surfaces$\}_{/\sim}$ is
the class of quasipositive surfaces up to ambient isotopy.

This paper is organized as follows.
In Section~2, we introduce the notion of front projections
in backslash position and prove that
a Legendrian ribbon is a quasipositive surface.
In Section~3, we give an alternative proof of Hedden's proposition.
Section~4 is devoted to Theorem~\ref{thmiso2} 
and Section~5 is devoted to Theorem~\ref{thmiso}.
In Section~6, we introduce the Thurston-Bennequin invariant
and the rotation number of a fence diagram of a quasipositive annulus
and prove their invariance under the moves of fence diagrams.

This work was done while the second author visited
Forschungsinstitut f\"ur Mathematik of ETH Z\"urich
and Mathematisches Institut of Universit\"at Basel
during the summer semester in 2006.
He would like to thank the staff of these departments,
especially Norbert A'Campo, for their warm hospitality
and precious support.

\section{Legendrian graphs and quasipositive diagrams}

The standard contact structure $\xi_{st}$ on $\R^3$ is the kernel 
of the $1$-form $dz+xdy$.
A {\it Legendrian graph} is a finite graph consisting
of edges and vertices such that
each edge is Legendrian, i.e.,
tangent to the $2$-plane field $\xi_{st}$ everywhere.
The graph may have a simple closed curve component,
which we also call an edge for convenience.
The image of the projection of a Legendrian graph $\Gamma$ onto 
the $yz$-plane is called a {\it front projection} of $\Gamma$.
This image is an immersed graph
with cusps and without vertical tangencies.
We call the image of each vertex also a vertex.
If $\Gamma$ is in general position,
its front projection has only node and cusp singularities,
and the edges adjacent to each vertex have the same tangency.
We call such a $\Gamma$ a {\it generic front projection}.
For each node, we regard the arc with smaller slope as
the strand passing over the other arc. Then the figure obtained
becomes a graph diagram of $\Gamma$.

A {\it Legendrian ribbon} $R$ of a Legendrian graph $\Gamma$
in $(\R^3,\xi_{st})$ is
a smoothly embedded surface in $S^3$ such that
\begin{itemize}
   \item[(1)]
$\Gamma$ is in the interior of $R$ and $R$ retracts onto $\Gamma$,
   \item[(2)]
for each $x\in \Gamma$, the $2$-plane of $\xi_{st}$ at $x$
is tangent to $R$, and
   \item[(3)]
for each $x\not\in \Gamma$,
the $2$-plane of $\xi_{st}$ at $x$ is transverse to $R$.
\end{itemize}
The notion of a ribbon of a Legendrian graph appears
in~\cite{giroux} to prove the existence of an open book
decomposition compatible with a given contact structure
on a $3$-manifold, see~\cite{etnyre,goodman}.

A quasipositive diagram retracts onto the preimage of its cusped
fence diagram. We call this preimage the {\it Legendrian core graph}.

\begin{lemma}\label{lemma2}
A quasipositive diagram can be regarded as a Legendrian ribbon 
of its Legendrian core graph.
\end{lemma}

\begin{figure}[htbp]
   \centerline{\input{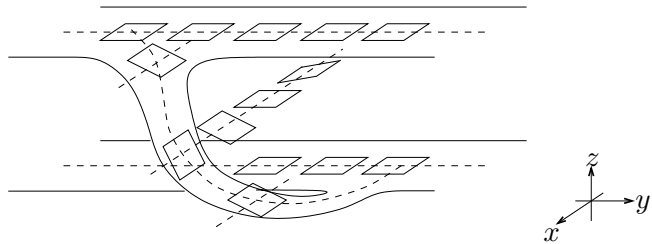}}
   \caption{A positive band of a quasipositive diagram in a position of
a Legendrian ribbon.\label{figure24}}
\end{figure}

\begin{proof}
Set the $b$ copies of disjoint parallel disks in $\R^3$ parallel
to the $xy$-plane and attach the positive bands in the region $x>0$
as shown in Figure~\ref{figure24}.
This figure shows that this surface is a Legendrian ribbon
of the Legendrian core graph.
\end{proof}

The main aim of this section is to proof the converse
of Lemma~\ref{lemma2}.

\begin{thm}\label{prop1}
A Legendrian ribbon of a Legendrian graph $\Gamma$ in $(\R^3,\xi_{st})$
is a quasipositive surface.
\end{thm}

Before proving the assertion, we introduce the notion of 
backslash position of front projections, trivalent front projections
and their approximating fence diagrams.

\begin{dfn}
A front projection is called {\it in backslash position}
if all tangent lines lie in $(\pi/2,\pi)\cup(3\pi/2,2\pi)$.
\end{dfn}

Define the diffeomorphism $\phi_1$ from the $yz$-plane
to itself by $(y,z)\mapsto (y, \lambda z)$,
where $\lambda>0$ is a positive real number,
and the diffeomorphism $\phi_2$ as the $-\pi/4$ rotation map
of the $yz$-plane.
For a given, Legendrian isotopy move of a front projection,
we can choose $\lambda$ sufficiently small such that
the diffeomorphism $\phi_2\circ\phi_1$
maps all front projections during the move into backslash position.

\begin{dfn}
A vertex in a generic front projection is called {\it trivalent}
if the number of adjacent edges is three
and two of the edges lie on one side with respect to the vertical line
passing through the vertex and the third edge lies on the other side.
A front projection is called {\it trivalent}
if it is generic and all vertices are trivalent.
\end{dfn}

\begin{figure}[htbp]
   \centerline{\input{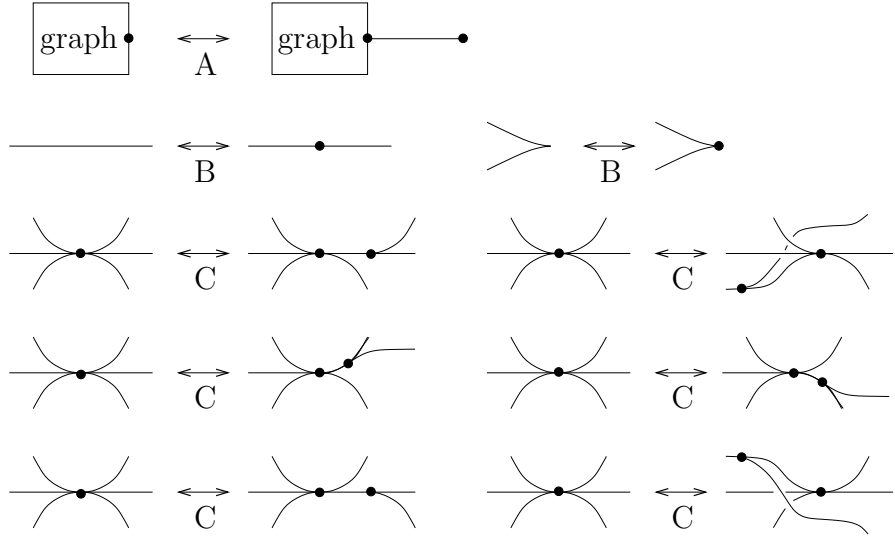}}
   \caption{The moves~A, B and~C of generic
front projections of Legendrian graphs.
The move B is allowed for left cusps.
In the figure of move C,
a slide of an edge on the left side of the vertex is also allowed.
\label{figure20}}
\end{figure}

We define the moves A, B and C of Legendrian graphs
as shown in Figure~\ref{figure20}.
The right direction of the move A is 
an addition of a pair of vertex $v$ and edge $e$ to a Legendrian graph.
The right direction of the move B is an addition of a vertex
to either the middle of an edge or to a cusp.
The move C is a slide of an edge.
The figure shows the case where the number of adjacent edges is six.

\begin{lemma}\label{lemma3}
The moves A, B and C satisfy the following properties:
\begin{itemize}
   \item[(1)]
the Legendrian ribbons before and after these moves
are ambient isotopic;
   \item[(2)]
every generic front projection can be modified into
a trivalent front projection by using these moves.
\end{itemize}
\end{lemma}

\begin{figure}[htbp]
   \centerline{\input{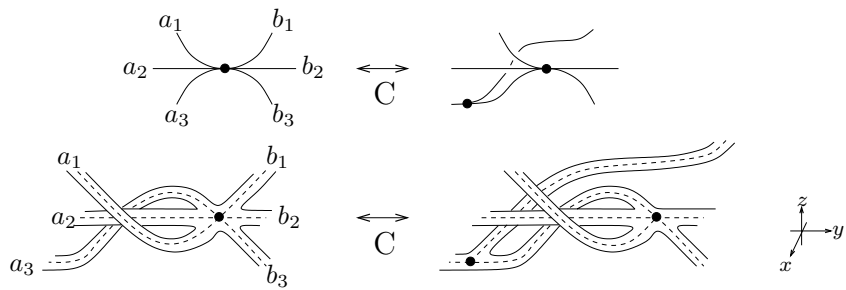}}
   \caption{The Legendrian ribbons before and after
the move C on the right-top in Figure~\ref{figure20}.\label{figure23}}
\end{figure}

\begin{proof}
The assertion (1) can be verified by describing their Legendrian ribbons.
The Legendrian ribbons before and after the move C on the right-top in
Figure~\ref{figure20} is shown in Figure~\ref{figure23}.
The assertion (2) is obvious.
\end{proof}

\begin{figure}[htbp]
   \centerline{\input{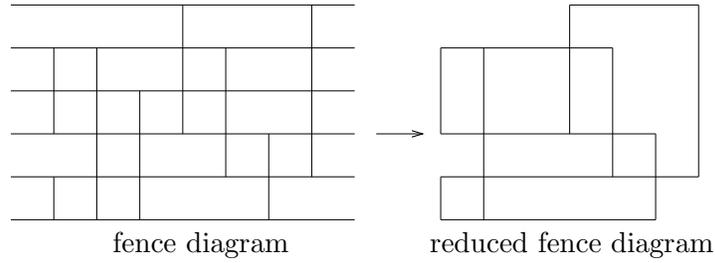}}
   \caption{A reduced fence diagram.\label{figure21}}
\end{figure}

\begin{dfn}
For a fence diagram, we apply deflations as much as possible
and then retract each of the left and right ends of horizontal lines
until their arriving at a trivalent vertex.
We call the obtained diagram the {\it reduced fence diagram}.
See Figure~\ref{figure21}.
The same operation is also applied to a cusped fence diagram
and we call the obtained diagram the {\it reduced, cusped fence diagram}.
\end{dfn}

Now we consider to approximate a trivalent
front projection in backslash position by reduced fence diagrams.
Let $q$ be a fence diagram, $r$ its reduced fence diagram
and $w$ a trivalent front projection in backlash position.
We denote by $\Sigma_c(r)$ the set of left-top and right-bottom
corners of $r$, which correspond to the cusps in the reduced, cusped
fence diagram of $q$, by $\Sigma_c(w)$ the set of cusps of $w$,
by $\Sigma_n(r)$ and $\Sigma_n(w)$
the set of nodes of $r$ and $w$ respectively,
and
by $\Sigma_v(r)$ and $\Sigma_v(w)$
the set of trivalent vertices of $r$ and $w$ respectively.

\begin{figure}[htbp]
   \centerline{\input{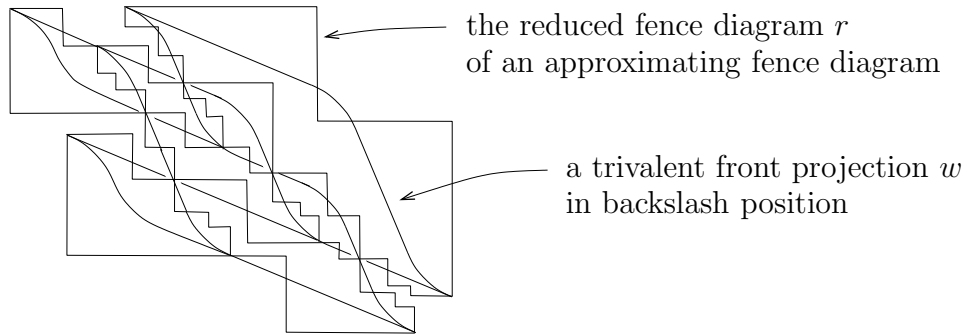}}
   \caption{An example of an approximating fence diagram.\label{figure16}}
\end{figure}

\begin{dfn}
We say a fence diagram $q$ {\it approximates} 
a trivalent front projection $w$ if its reduced fence diagram
$r$ satisfies the following:
\begin{itemize}
   \item[(1)] $\Sigma_c(r)=\Sigma_c(w)$;
   \item[(2)] $\Sigma_n(r)=\Sigma_n(w)$;
   \item[(3)] $\Sigma_v(r)=\Sigma_v(w)$;
   \item[(4)] there is a continuous family $r_t$ of curves
from $r_0=r$ to $r_1=w$, which is polygonal for $t=0$, such that
$r_{t_0}\cap r_{t_1}=\Sigma_c(r)\cup \Sigma_n(r)\cup (r\cap w)$
for all $0\leq t_0<t_1\leq 1$.
\end{itemize}
See Figure~\ref{figure16} for example.
\end{dfn}

\begin{lemma}\label{lemma28}
Let $q$ be a fence diagram, $r$ its reduced fence diagram
and $w$ a trivalent front projection in backslash position.
Suppose that $q$ approximates $w$. Then
\begin{itemize}
\item[(1)]
$r$ is regular isotopic to $w$ as immersed curves in $\R^2$ with cusps, and
\item[(2)]
the quasipositive surface of the fence diagram $q$
is a Legendrian ribbon of the Legendrian graph of $w$.
\end{itemize}
\end{lemma}

\begin{proof}
It is easy to verify the assertion (1), cf. Figure~\ref{figure16}. 
Isotope the quasipositive surface of the fence diagram $q$ 
according to the deflations and retractions for making the reduced
fence diagram and then isotope it further according to 
the isotopy move in the assertion~(1).
We then have a Legendrian ribbon of 
the Legendrian graph of $w$. This proves the assertion in (2).
\end{proof}

\vspace{3mm}

\noindent
{\it Proof of Theorem~\ref{prop1}.}\;\,
Let $\bar w$ be a generic front projection of a given Legendrian graph
and assume that $\bar w$ is in backslash position.
By Lemma~\ref{lemma3}, there exists a trivalent
front projection $w$ such that the Legendrian ribbons of $w$
and $\bar w$ are ambient isotopic.
By Lemma~\ref{lemma28},
the Legendrian ribbon of $w$ is quasipositive and
hence that of $\bar w$ is also.
\qed

\vspace{3mm}

\section{An alternative proof of Hedden's proposition}

Let $\alpha$ be a $1$-form on $S^3$ and $\xi=\ker\alpha$
its contact structure.
Two manifolds with contact structures are called {\it contactomorphic}
if there exists a diffeomorphism between these manifolds
which maps the $2$-plane field of the contact structure from
one to the other.
If $\xi$ is contactomorphic to the contact structure on
$S^3=\{(x_1,y_1,x_2,y_2)\in\R^4\mid x_1^2+y_1^2+x_2^2+y_2^2=1\}$
determined by the kernel of 
the $1$-form $\alpha=\sum_{i=1,2}(x_idy_i-y_idx_i)|_{S^3}$,
then $\xi$ is called the {\it standard contact structure} on $S^3$.
In this case, $(S^3,\xi)$ minus one point is contactomorphic
to $(\R^3,\xi_{st})$.

Let $F$ be an oriented surface with boundary
and $\phi:F\to F$ a diffeomorphism which is the identity map 
on the boundary $\bd F$ of $F$. 
Identify $F\times [0,1]$ by equivalence relations
$(x,1)\sim (\phi(x),0)$ for $x\in F$ and 
$(y,0)\sim (y,\theta)$ for $y\in \bd F$ and $\theta\in [0,1]$.
Then $F$ is called a {\it fiber surface}
in the $3$-manifold $F\times [0,1]/\sim$.
We consider only the case where $F\times [0,1]/\sim$ is $S^3$.
We denote the fiber surface parametrized by $\theta\in [0,1]$
by $F_\theta$ and suppose $F=F_0$.

A fiber surface $F$ embedded in $S^3$ 
is called {\it compatible} with a contact
structure $\xi$ on $S^3$ if it satisfies the following:
\begin{itemize}
   \item[(1)]
the boundary of $F_\theta$ is transverse to $\xi$;
   \item[(2)]
$d\alpha$ is a volume form on each fiber $F_\theta$;
   \item[(3)]
the orientation of $\bd F_\theta$ coincides with
the orientation of $\xi$ determined by $d\alpha>0$.
\end{itemize}

\begin{thm}[Hedden~\cite{hedden}]\label{thmhedden}
Let $F$ be a fiber surface in $S^3$.
$F$ is quasipositive if and only if $F$ is compatible with
the standard contact structure on $S^3$.
\end{thm}

We here give a proof of this theorem without using
the affirmative answer to Harer's conjecture.

\begin{proof}
Let $F$ be a fiber surface compatible with the standard contact structure
on $S^3$. By the Legendrian realization argument based
on the result in~\cite{giroux2}, we can assume that
the fiber surface is a Legendrian ribbon of a Legendrian graph,
cf.~\cite[Remark~4.30]{etnyre}.
Hence by Theorem~\ref{prop1} it is a quasipositive surface.

The proof of the converse assertion is done 
according to the argument in~\cite[p.21--23]{etnyre}.
Let $F$ be a quasipositive surface with a quasipositive diagram.
By Lemma~\ref{lemma2},
$F$ can be embedded in $S^3$ with the standard contact structure $\xi_0$
in such a way that $F$ is a Legendrian ribbon in $(S^3,\xi_0)$.
Since $\bd F$ is transverse to $\xi_0$,
there exists a small tubular neighborhood $N(\bd F)$ of $\bd F$ in $S^3$
with the contact structure $\ker(dz+r^2d\theta)$,
where the $z$-coordinate is along $\bd F$ and $(r,\theta)$ is
the polar coordinates of a plane transverse to $\bd F$.
We then define $M$ to be the union of $F\times[-\ve,\ve]$ and $N(\bd F)$,
where $\ve$ is a sufficiently small positive real number,
and assume that the boundary of $M$ is convex
(see~\cite{giroux2}, or for instance~\cite{eh} for
the definition of convexity).
Since $F$ is a fiber surface, 
the complement $M^c=\text{closure}(S^3\setminus M)$
is a handlebody with the tight contact structure $\xi_0|_{M^c}$.

The rest of the proof is the same as
the argument in~\cite[p.21--23]{etnyre}, so we only show the outline.
We first deform the Reeb vector field of $(S^3,\xi_0)$ in such a way that
it is tangent to the boundary of $N(\bd F)$ and transverse
to the fibers of the fibration in $F\times[-\ve,\ve]$.
Next we once forget the contact structure $\xi_0$ on $M^c$
and extend the Reeb vector field on $M$ to $M^c$ according to the 
fibration. In particular, the contact structure $\xi_1$
determined by this Reeb vector field is compatible with the fibration.
The Reeb vector field allows us to make a contact embedding
of $(M^c,\xi_1)$ into
$F\times\R$ with the vertically invariant contact structure.
By Giroux's criterion, the contact structure on $F\times\R$ is tight
and hence $(M^c,\xi_1)$ is also.
Since two contact structures $\xi_0|_{M^c}$ and $\xi_1$ on $M^c$
are both tight, due to the uniqueness of the tight contact structure
on a handlebody~\cite{torisu},
we can conclude that $\xi_0|_{M^c}$ and $\xi_1$ are contactomorphic.
This means that the contact structure $(M,\xi_0|_M)\cup(M^c,\xi_1)$
is contactomorphic to $(S^3,\xi_0)$, which is the standard
contact structure on $S^3$.
\end{proof}

\section{Isotopy moves of Legendrian graphs and quasipositive diagrams}

\begin{figure}[htbp]
   \centerline{\input{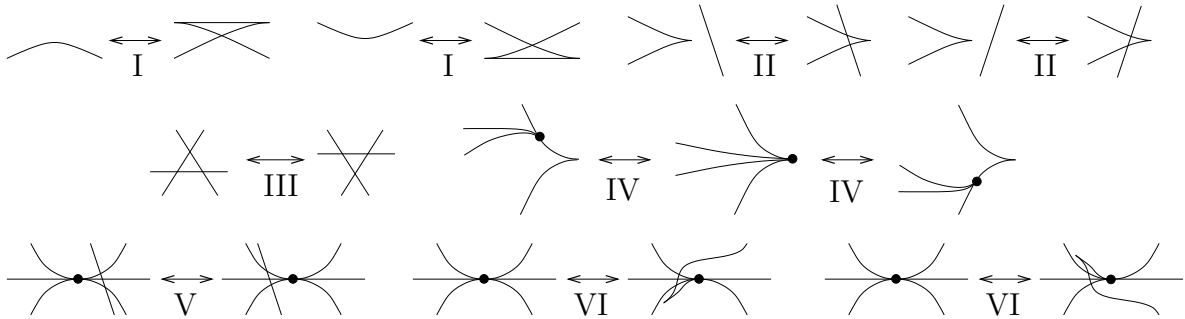}}
   \caption{Legendrian isotopy moves. The horizontal reflections
of these moves are also allowed. The overstrand and understrand
at each crossing are determined according to the rule that
the arc with the smaller slope passes over the other arc.\label{figure12}}
\end{figure}

It is well-known that 
two generic front projections of a Legendrian knot are related by
the moves I, II and~III of generic front projections shown
in Figure~\ref{figure12} (\cite{swiatkowski}).
In case of Legendrian graphs, we can assume that
the vertices are in general position during Legendrian isotopy moves
so that they do not intersect.
Mutual positions of a vertex and edges during a Legendrian isotopy move
yield three additional cases:
(IV)~a vertex passes through a cusp,
(V)~a vertex passes over or under an edge, and
(VI)~an edge adjacent to a vertex rotates to the other side of the vertex.
These moves are also described in Figure~\ref{figure12}.

\begin{figure}[htbp]
   \centerline{\input{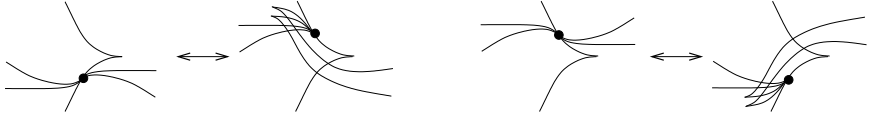}}
   \caption{The move~IV with edges attached to the right side
of the cusp.\label{figure28}}
\end{figure}

\begin{rem}\label{rem28}
The move~IV with edges attached to the right side of the cusp
as shown in Figure~\ref{figure28} is realized as a combination
of the moves~VI and~IV.
\end{rem}

\begin{thm}\label{thmiso2}
If the Legendrian graphs of reduced, cusped fence diagrams
are Legendrian isotopic, then
their fence diagrams are related by inflations,
deflations, slips and slides.
\end{thm}

\begin{rem}\label{remdfn}
In~\cite{rudolph:98}, the slides were defined as the moves
from the left to the right in Figure~\ref{figure3}
and it was remarked on p.263 that
the inverse moves can be realized as conjugations of slides
by twirls. In this paper, for our convenience to
compare with Legendrian isotopy moves,
we call these inverse moves also slides.
The definitions of twirls and turns in Figure~\ref{figure3}
are also different from those in~\cite{rudolph:98}
because of the same reason.
\end{rem}

\begin{figure}[htbp]
   \centerline{\input{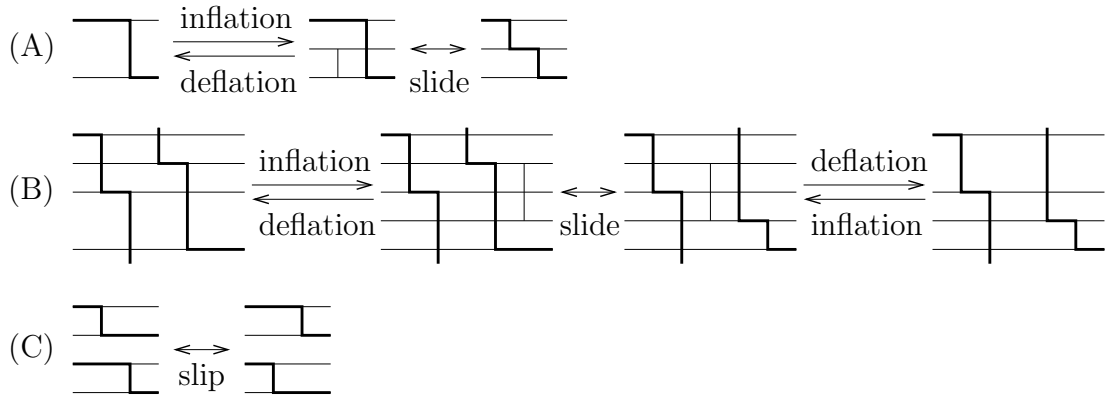}}
   \caption{A new zigzag and exchanges of horizontal and vertical edges.\label{figure17}}
\end{figure}

Before proving Theorem~\ref{thmiso2},
we explain the flexibility of approximating fence diagrams
shown in Figure~\ref{figure17}.
The thickened polygonal curves in the figures represent a part of
the reduced fence diagram. Figure (A) shows that
by combining an inflation and a slide
we can produce a new zigzag of an approximating fence diagram.
Figure (B) shows that by combining an inflation, slides and a deflation
we can exchange the heights of two horizontal edges,
and Figure (C) shows that
by a slip we can exchange the positions of vertical edges.
These properties imply that
we can make an approximating fence diagram
whose reduced fence diagram is as close to $w$ as
we need
by using inflations and slides, and
every regular isotopy move of a generic front projection,
as moves of immersed curves in $\R^2$ with cusps,
can be expressed by a family of approximating fence diagrams
defined by inflations, deflations, slips and slides.

\vspace{3mm}

\noindent
{\it Proof of Theorem~\ref{thmiso2}.}\;\,
We first show that the Legendrian isotopy moves I $\sim$ VI
can also be expressed by moves of approximating fence diagrams.
Since all reduced fence diagrams are trivalent,
the vertices in the moves IV, V and VI are trivalent.
Consider the move~II with right cusp.
If the cusp passes a line from the top to the bottom, then the move
is expressed by a combination of an inflation and slides
as shown on the top in Figure~\ref{figure18}.
\begin{figure}[htbp]
   \centerline{\input{figure18.pstex_t}}
   \caption{The moves of approximating fence diagrams corresponding
to the Legendrian isotopy moves II, III, IV, V and VI.\label{figure18}}
\end{figure}
The move of approximating fence diagrams corresponding to
the move~II with left cusp is given by the $\pi$-rotation of the figure.
In case the right (resp. left) cusp passes a line from
the left to the right (resp. from the right to the left), the move
corresponds to a slip, cf. Figure~\ref{figure15} below.
For every move which appears below, we also need to check the figure
obtained by the $\pi$-rotation, though we omit it.
The move~III corresponds to a slip, which is shown 
on the second figure in Figure~\ref{figure18} and
the move~IV corresponds to a slide as shown on the third figure.
The move V also corresponds to a slide if
the cusp passes a line from the top to the bottom, see the fourth figure.
In case where the cusp passes a line from the left to the right, the move
corresponds to a slip as in case of the move~II.
The move VI has four cases as shown in the fifth figure.
The move of approximating fence diagrams for the first case
is shown on the bottom in Figure~\ref{figure18}.
Looking only at the cusp in the move VI on the bottom
in Figure~\ref{figure18}, we have the move I.
We can check the other cases of the move~VI by the same way.
Thus we conclude that the moves~I $\sim$ VI are expressed
by approximating fence diagrams with using
inflations, deflations, slips and slides.

Now we prove the assertion.
Let $w$ and $w'$ be reduced, cusped fence diagrams
whose Legendrian graphs are Legendrian isotopic and
$r_0$ and $r'_0$ 
the reduced fence diagrams corresponding to $w$ and $w'$ respectively.
By small perturbation we can assume that
$w$ and $w'$ are in backslash position, and
by adding new zigzags to $r_0$ and $r'_0$
we can obtain fence diagrams $r_1$ and $r'_1$
which approximate $w$ and $w'$ respectively.
Now we make the approximation $r_1$ of $w$ as close to $w$ as possible
and follow the Legendrian isotopy moves from $w$ to $w'$
with approximating fence diagrams.
We denote the obtained fence diagram by $r_2$, which approximates $w'$.
It is easy to make the same approximation of $w'$ from
$r'_1$ and $r_2$ by adding new zigzags.
Thus $r_0$ and $r'_0$ are related by the moves in the assertion.
\qed

\vspace{3mm}

\section{Quasipositive annuli}

In this section we study quasipositive annuli.
The reduced fence diagram of a quasipositive diagram of a
quasipositive annulus has no trivalent vertices,
i.e., it is a knot diagram of a knot in $\R^3$.

\begin{thm}\label{thmiso}
The moves of fence diagrams of quasipositive annuli
under inflation, deflation, slips and slides
correspond to Legendrian isotopy moves of reduced, cusped fence diagrams. 
\end{thm}

\begin{proof}
An inflation and a deflation do not change the reduced
fence diagram. The slip of a reduced fence diagram
shown in Figure~\ref{figure15} corresponds to the Legendrian move II
if there is no horizontal line passing under the shorter vertical edge.
In case where the shorter vertical edge passes over several
horizontal lines, the move is realized by the moves II and III.
The other cases of slips are obviously Legendrian isotopy.

\begin{figure}[htbp]
   \centerline{\input{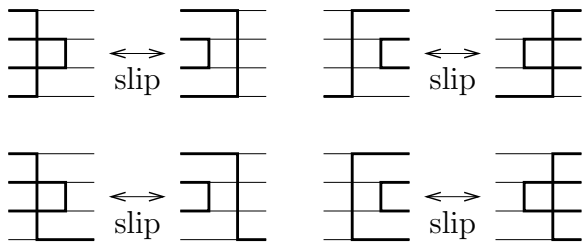}}
   \caption{Moves of reduced fence diagrams under slips.\label{figure15}}
\end{figure}

We consider the slide from the left to the right 
shown in Figure~\ref{figure5} (A).
If the reduced fence diagram does not
pass through the lower vertical edge in Figure~\ref{figure5}~(B),
then this move is obviously Legendrian isotopy.
If it passes through the lower vertical edge, there are eight cases
shown in Figure~\ref{figure5}~(C).
The non-obvious case is (C3) and corresponds
to the Legendrian isotopy move~I.
In case where the vertical lines of fence diagrams in Figure~\ref{figure5}
pass over several other horizontal lines,
we need to use the Legendrian isotopy move II additionally.
The proof for the slide from the right to the left also
follows from Figure~\ref{figure5}.
For the second slide in Figure~\ref{figure3}, the proof is analogous.
\end{proof}

\begin{figure}[htbp]
   \centerline{\input{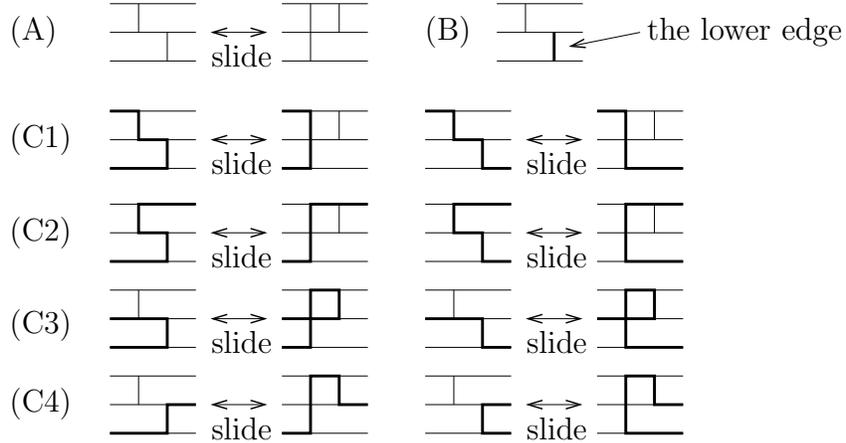}}
   \caption{Moves of reduced fence diagrams under slides.\label{figure5}}
\end{figure}

\begin{figure}[htbp]
   \centerline{\input{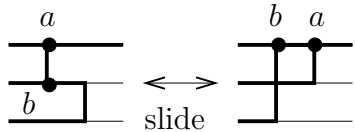}}
   \caption{A slide which exchanges the mutual
positions of two vertices.\label{figure27}}
\end{figure}

\begin{rem}
For quasipositive surfaces,
the moves of their fence diagrams under inflation, deflation and slips
correspond to Legendrian isotopy moves of reduced, cusped fence diagrams. 
But this is not true for slides. Figure~\ref{figure27} shows that
a slide may exchange the mutual positions of two vertices
and this cannot be realized by Legendrian isotopy moves
of Legendrian graphs.
\end{rem}

\begin{rem}
Theorem~\ref{thmiso2} shows the existence of the surjective map
\[
\{\text{trivalent Legendrian ribbons}\}_{/\sim}
\to
\{\text{quasipositive diagrams}\}_{/\sim},
\]
where 
$\{$trivalent Legendrian ribbons$\}_{/\sim}$
is the class of Legendrian graphs up to Legendrian isotopy
and
$\{$quasipositive diagrams$\}_{/\sim}$ is
the class of quasipositive diagrams up to
inflations, deflations, slips, slides, twirls and turns,
and the above remark shows that this map is not injective.
\end{rem}

\section{Quasipositive surfaces with different fence diagrams}

Let $\gamma$ be a Legendrian knot in $(\R^3,\xi_{st})$, i.e.,
$\gamma$ is tangent to the $2$-plane field $\xi_{st}$ everywhere.
The {\it Thurston-Bennequin invariant} ${\it tb}(\gamma)$ of 
$\gamma$ is the linking number of $\gamma$
and a curve obtained by pushing off $\gamma$ normal to $\xi_{st}$.
To give the definition of the rotation number,
we assign an orientation to $\gamma$ and denote it by $\vec\gamma$.
The {\it rotation number} $r(\vec\gamma)$ of $\vec\gamma$ is
the winding number of vectors tangent to $\vec\gamma$
with respect to the trivialization of $\xi_{st}$ along $\vec\gamma$.
The other choice of the orientation of $\gamma$ changes
the sign of the rotation number.

The Thurston-Bennequin invariant and the rotation number
can be read from the front projection.
Let $\gamma$ be a Legendrian knot
and $w_\gamma$ its generic front projection.
Assign an orientation to $w_\gamma$ and
let $p(w_\gamma)$ (resp. $n(w_\gamma)$) denote 
the number of positive (resp. negative) crossings
and $r_c(w_\gamma)$ the number of right cusps of $w_\gamma$.
Then the Thurston-Bennequin invariant of $\gamma$
is determined by the formula
\begin{equation}\label{tb}
   {\it tb}(\gamma)=p(w_\gamma)-n(w_\gamma)-r_c(w_\gamma).
\end{equation}
Obviously, this number does not depend on the choice of the orientation.
For the rotation number, 
let $d_c(w_\gamma)$ (resp. $u_c(w_\gamma)$) denote
the number of downward (resp. upward) cusps.
Then the rotation number is determined by
\[
   r(\vec\gamma)=\frac{1}{2}(d_c(w_\gamma)-u_c(w_\gamma)).
\]
It is easy to verify that the other choice of the orientation of
$\gamma$ changes the sign of the rotation number.
For more precise explanations, see for instance~\cite{gompf:98}.

A reduced, cusped fence diagram of a quasipositive annulus
is regarded as a front projection of a Legendrian knot.
We define the Thurston Bennequin invariant
and the rotation number of a fence diagram of a quasipositive annulus
by those of its reduced, cusped fence diagram.

\begin{lemma}\label{lemma1}
The Thurston-Bennequin invariant of a fence diagram of a
quasipositive annulus $A$ is equal to $-1$ times the linking number of
the two boundary components of $A$.
In particular, the Thurston-Bennequin invariant is independent
of the choice of a quasipositive diagram.
\end{lemma}

\begin{figure}[htbp]
   \centerline{\input{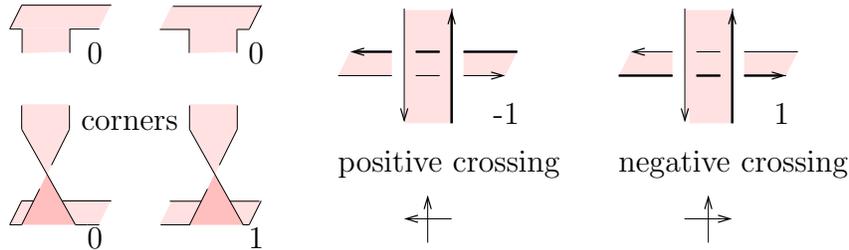}}
   \caption{The contribution to the linking number of
the two boundary components of $A$.\label{figure14}}
\end{figure}

\begin{proof}
On a quasipositive diagram, we assume that the positive crossing
of each positive band lies close to the bottom end of the band.
Then the contribution to the linking number of the two boundary
components of $A$ is given as shown in Figure~\ref{figure14}.
The number on the right-bottom of each figure represents
the contribution to the linking number.
Thus the linking number is
$-p+n+r_c=-tb$ by the formula~\eqref{tb}.
\end{proof}

\begin{rem}
Lemma~\ref{lemma1} can be proved by checking the coincidence of
the Legendrian framing and the Seifert surface's framing of
a reduced, cusped fence diagram,
cf. for example~\cite{ao:2001}.
This proof is more direct than the above proof
if we assume the knowledge of these framings.
\end{rem}

\begin{thm}\label{mainthm}
The Thurston-Bennequin invariant and the rotation number
of a fence diagram of a quasipositive annulus
are invariant under inflations, deflations, slips, slides, twirls and turns.
\end{thm}

\begin{proof}
Since the moves of quasipositive annuli in the assertion
are ambient isotopy moves of Seifert surfaces,
the linking number of the two boundary components of 
a quasipositive annulus does not change under these moves.
Therefore, by Lemma~\ref{lemma1},
the Thurston-Bennequin invariant also does not change.
The rotation number is also invariant
under inflations, deflations, slips and slides
since they are Legendrian isotopy moves by Theorem~\ref{thmiso}.



We will prove the invariance of the rotation number
in case of twirls and turns.
For twirls,
there are four cases (A), (B), (C) and (D) of reduced fence diagrams
as shown in Figure~\ref{figure7}.
We assign an orientation as shown in the figures.
The small arrows in the figures represent the positions
of upward and downward cusps. It is easy to check that
the rotation number does not change under these moves.
The other choice of the orientation changes the sign of
the rotation number, but it does not matter for
its invariance under these moves.

\begin{figure}[htbp]
   \centerline{\input{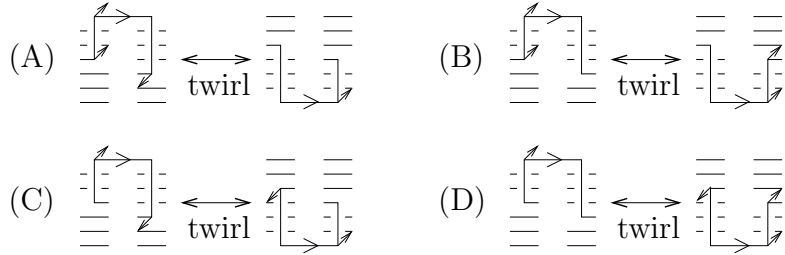}}
   \caption{The four cases of the moves of reduced fence diagrams under twirls.\label{figure7}}
\end{figure}

The proof of the invariance under turns is
analogous to the proof for twirls.
There are four cases (A), (B), (C) and (D) of the moves
of reduced fence diagrams as shown in Figure~\ref{figure8}
and we can check easily that the rotation number does not change
under these moves.
\end{proof}

\begin{figure}[htbp]
   \centerline{\input{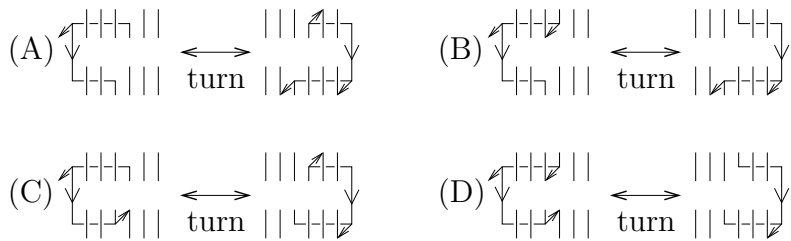}}
   \caption{The four cases of the moves of reduced fence diagrams under turns.\label{figure8}}
\end{figure}

As a corollary, we answer a question of Rudolph
in~\cite[Remark in p.263]{rudolph:98}.

\begin{cor}\label{maincor}
There exists a quasipositive surface with
two different quasipositive diagrams which are not related by
inflations, deflations, slips, slides, twirls and turns.
\end{cor}

\begin{figure}[htbp]
   \centerline{\input{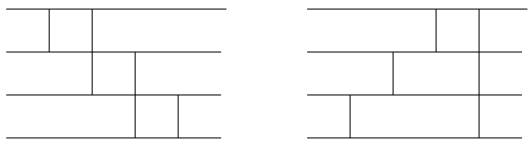}}
   \caption{Two different fence diagrams of the $3$ times full twisted quasipositive annulus.\label{figure10}}
\end{figure}

\begin{proof}
We consider two fence diagrams shown in Figure~\ref{figure10}.
They are quasipositive annuli and
it is easy to check by the formula~\eqref{tb} that
the Thurston-Bennequin invariants are both $-3$.
Hence, by Lemma~\ref{lemma1},
their quasipositive surfaces are the same, namely
the $3$ times full twisted quasipositive annulus.
However the rotation number of the fence diagram on the left
is $0$ and the number of the right diagram is $\pm 2$. 
Hence, by Theorem~\ref{mainthm},
these fence diagrams are not related by the moves in the assertion.
\end{proof}

If the Legendrian isotopy class of a Legendrian knot
with the same knot-type is determined by the Thurston-Bennequin invariant
and the rotation number, then this knot-type 
is called {\it Legendrian simple}.
It is known by Y.~Eliashberg and M.~Fraser in~\cite{ef}
that the unknot is Legendrian simple, and
J.~Etnyre and K.~Honda proved in~\cite{eh} that
the torus knots and the figure eight knot are Legendrian simple.
On the other hand, the knots $5_2$, $6_3$ and $7_2$,
in Rolfsen's notation~\cite{rolfsen}, are not
Legendrian simple~\cite{chekanov,ng}. These facts are known 
as an application to Chekanov's differential graded algebra~\cite{chekanov}.

As a direct corollary of Theorem~\ref{mainthm},
we can determine the isotopy class of a quasipositive annulus
up to the moves of quasipositive annuli
in case where its core curve is Legendrian simple.

\begin{cor}\label{cor3}
Let $A$ be a quasipositive annulus such that
the knot type of its core curve is Legendrian simple.
\begin{itemize}
   \item[(1)]
The isotopy class of a quasipositive diagram of $A$
up to inflations, deflations, slips and slides
is determined by the Thurston-Bennequin invariant and
the rotation number of their reduced, cusped fence diagrams.
   \item[(2)]
The classification of quasipositive diagrams of $A$ in (1) 
is equivalent to the classification 
up to inflations, deflations, slips, slides, twirls and turns.
\end{itemize}
\end{cor}

For example,
consider the $n$ times full twisted quasipositive annulus $A_n$.
By using the classification of the Legendrian unknot in~\cite{ef},
we know that there exist $\lfloor (n+1)/2 \rfloor$ reduced
fence diagrams of $A_n$ with different rotation numbers.
Hence, by Corollary~\ref{cor3},
we conclude that they are not related by
inflations, deflations, slips, slides, twirls and turns.

\begin{rem}
The existence of the map
\[
\{\text{quasipositive diagrams}\}_{/\sim}
\to
\{\text{quasipositive surfaces}\}_{/\sim}
\]
is obvious and Corollary~\ref{maincor} shows that this map is not injective,
where $\{$quasipositive surfaces$\}_{/\sim}$ is
the class of quasipositive surfaces up to ambient isotopy and
$\{$quasipositive diagrams$\}_{/\sim}$ is
the class of quasipositive diagrams up to
inflations, deflations, slips, slides, twirls and turns.
\end{rem}

\begin{q}
Is there a quasipositive fiber surface
with different fence diagrams up to 
inflations, deflations, slips, slides, twirls and turns?
\end{q}

\begin{q}
Do the moves of reduced fence diagrams under 
twirls and turns correspond to Legendrian isotopy moves?
\end{q}

\begin{q}
Follow the proof of Harer's conjecture by using fence diagrams.
\end{q}

\begin{q}
Does there exist a quasipositive fiber surface, other than a disk,
from which we cannot deplumb a Hopf band?
\end{q}

\end{document}